

\baselineskip=14pt
\parskip=10pt

\magnification=\magstephalf

\def\1{{\overline{1}}}
\def\2{{\overline{2}}}
\parindent=0pt
\overfullrule=0in

\def\frac#1#2{{#1 \over #2}}

\centerline
{
\bf The Binomial Theorem  for $(N+n)^r$ (where Nf(n)=f(n+1))
}
\bigskip
\centerline
{\it Moa APAGODU, Shalosh B. EKHAD, and Patrick GASKILL}
\bigskip

{\bf Added Dec. 13, 2011}: The main theorm of this note is contained in Lemma 5
of ``On Two-generated Non-commutative Algebras Subject to the Affine Relation''
by Christoph Koutschan, Viktor Levandovskyy, Oleksandr Motsak, {\tt http://arxiv.org/abs/1108.1108},
who prove many other results, and a stronger version of our main result
(using Stirling numbers).

This note has only personal and historical interest, and is only published in
the Personal Journal of Ekhad and Zeilberger and arxiv.org .
\bigskip
\hrule
\bigskip
We all know the {\it binomial theorem}
$$
(x+y)^r=\sum_{i=0}^{r} \frac{r!}{i!(r-i)!} x^i y^{r-i} \quad ,
\eqno(1)
$$
where $x$ and $y$ are {\it commuting} variables, i.e. $yx=xy$.
The binomial theorem is  easily proved by induction on $r$.

Not as famous is the {\it quantum} analog, that goes back to Marco Sch\"utzenberger,
$$
(x+y)^r=\sum_{i=0}^{r} \frac{[r]!}{[i]![r-i]!} x^i y^{r-i} \quad ,
\eqno(2)
$$
where $x$ and $y$ are {\it q-commuting} variables, i.e. $yx=qxy$,
and $[j]!:=(1)(1+q)(1+q+q^2) \cdots (1+q+ \dots +q^{j-1})$.
The quantum binomial theorem is also easily proved by induction on $r$.

Even more obscure is the binomial theorem for $(D+x)^r$, 
where $D$ is the differentiation operator $\frac{d}{dx}$
(so $Dx=xD+1$):
$$
(x+D)^r= \sum_{k=0}^{\lfloor r/2 \rfloor}
\frac{r!}{2^k k!}
\sum_{j+l=r-2k} \frac{1}{j!l!} x^j D^l \quad ,
\eqno(3)
$$
that is also easily proved by induction.

But nothing analogous is known for $(n+N)^r$, where $N$ is the {\it shift operator} $Nf(n):=f(n+1)$,
and $n$ is multiplication by $n$. Now the {\it commutation relation} is $Nn=nN+N$ and things
seem to get messier.

Let's first try to expand $(N+n)^r$ in powers of $N$ with coefficients that are polynomials in $n$:
$$
(N+n)^r=\sum_{d=0}^{r} P_{r,d}(n) N^d \quad .
\eqno(4a)
$$
We claim that the coefficients $P_{r,d}(n)$ are given by the following ``explicit'' formula
$$
P_{r,d}(n)=
\sum_
{
{
{p_0+\dots+p_d=r-d}
\atop
{0 \leq p_0, \dots p_d \leq r-d}
}
}
\, \prod_{j=0}^{d} (n+j)^{p_j} \quad .
\eqno(4b)
$$
In other words,  $P_{r,d}(n)$ is the {\it weight-enumerator} of {\it compositions} of $r-d$ into
$d+1$ {\it non-negative} integers, with 
$$
Weight([p_0, \dots, p_d]):=
\prod_{j=0}^{d} (n+j)^{p_j} \quad .
$$

Since $(N+n)^r=(N+n)(N+n)^{r-1}$ we have the recurrence:
$$
P_{r,d}(n)=P_{r-1,d-1}(n+1)+nP_{r-1,d}(n) \quad .
$$
Identity $(4)$ is proved by induction by noting that any composition of 
$r-d$ into $d+1$ non-negative integers either
has $p_0=0$ and beheading it yields a composition of $r-d$ into $(d-1)+1$ 
non-negative integers,
and the weights get adjusted by replacing $n$ by $n+1$,
or $p_0 \geq 1$ and 
subtracting $1$ from $p_0$ yields a 
composition of $r-d-1$ into $d+1$ non-negative integers,
and adding the $1$ back to the $p_0$ term results in multiplying the weight by $n$.

While formula $(4)$ is very elegant and combinatorial, 
it would be nice to have an {\it explicit}
formula, as a linear combination of monomials $n^iN^j$ analogous to 
$(1)$, $(2)$ and $(3)$.
This does not seem to be possible, but by using the Maple package  
{\tt Nnr}, written
by Doron Zeilberger, and downloadable from 
{\tt http://www.math.rutgers.edu/{}\~{}zeilberg/tokhniot/Nnr},
one can get an explicit formula for the first $k$ highest-degree 
terms for any desired $k$, for $r \geq 2k$. 

Let's describe the answer for $k=3$. First we define,
$$
((n+N))^r=\sum_{i=0}^{r} \frac{r!}{i!(r-i)!} n^i N^{r-i} \quad ,
$$
in other words the polynomial in the (non-commuting) variables $n$ and $N$ 
obtained
by expanding $(n+N)^m$  while pretending that $n$ and $N$ commute. We have 
(all the algebric computations below should be done in commutative algebra,
and at the end  each monomial should be expressed with $n$ before $N$. i.e. in the form $n^iN^j$)
$$
(N+n)^r=
((N+n))^r+
{{r} \choose {2}}N((N+n))^{r-2}
+
{{r} \choose {3}}N((N+n))^{r-4}
\left (\frac{1}{4} \, \left( 3\,r-5 \right) N+n \right )
$$
$$
+
{{r} \choose {4}}N ((N+n))^{r-6}
\left (
\frac{1}{2}\, \left( r-2 \right)  
\left( r-3 \right) {N}^{2}+2\, \left( r-3 \right) nN+{n}^{2}
\right )
\quad +
$$
$$
(terms-of-total-degree \,\,\, \leq r-4) \quad .
$$
To get all the terms of total degree $\geq r-10$, see 
{\tt http://www.math.rutgers.edu/\~{}zeilberg/tokhniot/oNnr10}.
You are welcome to get all the terms of degree 
$\geq r-k$, for any desired positive
integer $k$, by typing {\tt Mispat(r,k,n,N); } in the Maple package {\tt Nnr}.

\bigskip
\hrule
\bigskip

Moa Apagodu,
Department of Mathematics and Applied Mathematics,
Virginia Commonwealth University, Richmond, VA, 23284, 
{\tt mapagodu [at] vcu [dot] edu} .

\bigskip
Shalosh B. Ekhad, Department of Mathematics, 
Rutgers University, Piscataway, NJ 08854 .

\bigskip
Patrick Gaskill,
Department of Mathematics and Applied Mathematics,
Virginia Commonwealth University, Richmond, VA, 23284, 
{\tt gaskillpw [at] vcu [dot] edu} .

\end